\documentclass{amsart}

\usepackage[latin1]{inputenc}
\usepackage{latexsym}
\usepackage{amsfonts}
\usepackage{color}
\usepackage{dsfont}
\usepackage{varioref}
\usepackage{epsfig}
\usepackage{hyperref}
\usepackage[english]{babel}

\newcommand{\dP}{{\mathds P}}
\newcommand{\dC}{{\mathds C}}
\newcommand{\dR}{{\mathds R}}
\newcommand{\dQ}{{\mathds Q}}
\newcommand{\smcdot}{{\textup{$\cdot$}}}

\newtheorem{theorem}{Theorem}
\newtheorem{corollary}[theorem]{Corollary}

\begin{document}

\title{A Sextic with $35$ Cusps}
\author{Oliver Labs}
\address{Johannes Gutenberg Universit\"at Mainz,
  Germany}
\email{Labs@Mathematik.Uni-Mainz.de, mail@OliverLabs.net}
\date{\today}

\subjclass{Primary 14J17, 14Q10}
\keywords{algebraic geometry, many cusps, many singularities}


\begin{abstract}
  Recently, W.~Barth and S.~Rams discussed sextics with up to $30$ 
  $A_2$-singularities (also called cusps) and their connection to coding
  theory \cite{barramsCC}. 
  In the present paper, we find a sextic with $35$
  cusps within a four-parameter family of surfaces of degree $6$ in projective
  three-space with dihedral symmetry $D_5$.    
  This narrows the possibilities for the maximum number $\mu_{A_2}(6)$ 
  of $A_2$-singularities on a sextic to $35 \le \mu_{A_2}(6) \le 37$.
  To construct this surface, we use a general algorithm in characteristic zero
  for finding hypersurfaces with many singularities within a family. 
\end{abstract} 

\maketitle

\section*{Introduction}

Since the middle of the $19^\textup{\tiny th}$ century algebraic
geometers are interested in the question: 
How many isolated singularities of a given topological type can a surface of
given degree $d$ in $\dP^3 := \dP^3(\dC)$ contain?  
For surfaces of degree $d\le3$ the complete answer to this question is already
known since Schl\"afli's work \cite{Schl63} on cubic surfaces in 1863, see
\cite{labsHolzerCS} for explicit equations and illustrating pictures.
The treatment of the case $d=4$ was completed in 1997 using computers (see
\cite{yangListQuartics} and \cite{urabeInv}). 

For higher degree much less is known, even when restricting to surfaces with
$A_k$-singularities with local equation $x^{k+1}+y^2+z^2=0$.
We denote by $\mu_{A_k}(d)$ the maximum number of $A_k$-singularities a
surface of degree $d$ in $\dP^3$ can have. 
The maximum number $\mu_{A_1}(d)$ of nodes on a surface of degree $d$
is only known for $d\le6$ (see \cite{bar65} and \cite{jafrub66} for the case
$d=6$).  
In \cite{labs99} we improved the lower bound for the next open case by
constructing a surface of degree $7$ with $99$ nodes which narrowed the
possibilities for $\mu_{A_1}(7)$ to: $99\le\mu_{A_1}(7)\le104$.   
For most higher degrees, Chmutov \cite{chmuP3} described the currently best known
construction.  

Even less is known for $A_2$-singularities.
The currently best known upper bounds follow from Miyaoka's result
\cite{miyP3} and for $d=4,7$ from \cite{varBound}.    
The lower bounds are usually  achieved by constructions.
The cases $d\le3$ are classical.
$\mu_{A_2}(4)=8$ follows from Yang's article mentioned above in an abstract
way, and Barth gave an explicit construction for an $8$-cuspidal quartic
in \cite{barClassK3}. 
In \cite{barQuint}, he showed $\mu_{A_2}(5)\ge15$ by constructing a
$15$-cuspidal quintic which is nicely connected to Clebsch's Cubic Diagonal
Surface.  
For the construction of sextics with up to $30$ cusps (this was the maximum
known prior to the present article), based on an idea
already used by Rohn in the $19^\textup{\tiny th}$ century, see
\cite{barramsCC}.  
In this article, Barth and Rams also study the codes connected to these
cuspidal sextics.  

We also base our construction on Rohn's idea.  
This gives a $4$-parameter family of sextics with
dihedral symmetry $D_5$ and $30$ cusps. 
Then we use the computer algebra programm {\sc Singular}
\cite{Singular} to find a sextic with $35$ cusps within this family 
which shows: $$\mu_{A_2}(6) \ge 35.$$

The following table lists the known restrictions to $\mu_{A_2}(d)$. 
For $d\ge7$, they follow from our generalization of Chmutov's construction of
nodal surfaces to higher $A_j$-singularities, see \cite{labsVarChmu}. 

\begin{center}
\begin{tabular}{|p{0.9cm}|p{0.15cm}|p{0.15cm}|p{0.15cm}|p{0.35cm}|p{0.35cm}|p{0.35cm}|p{0.35cm}|p{0.45cm}|p{0.45cm}|p{0.45cm}|p{1.65cm}|}
\hline
\rule{0pt}{0.95em}degree & 2 & 3 & 4 & \ 5 & \ 6 & \ 7 & \ 8 & \ \ 9 & \ 10
& \ 11 & d\\
\hline
\hline
\rule{0pt}{1.05em}$\mu_{A_2} \le$ & $0$ & $3$ & $8$ & $20$ & $37$ & $62$ & $98$ & $144$ & $202$
& $275$ & ${1\over4}d(d-1)^2$\\[0.1em] 
\hline
\rule{0pt}{1.15em}$\mu_{A_2} \ge$ & $0$ & $3$ & $8$ & $15$ & {\bf 35} &
$52$ & $70$ & $126$ & $159$ & $225$ &
$\approx {2\over9}d^3$\\[0.1em] 
\hline
\end{tabular}
\end{center}

\vspace{0.6em}
I thank D.~van~Straten and W.~Barth for valuable discussions and W.~Barth for
his visit to Mainz which was a good motivation to complete this work.

\section{The Sextics with $35$ Cusps}
\label{secSextics}

In previous works (e.g., \cite{vStrQuint}, \cite{bar65},
\cite{endrOct}), the authors used geometric
arguments to reduce a problem depending on several parameters to 
polynomials each depending only on one parameter.  
The roots of these polynomials could then easily be found by hand or by
computer algebra. 
But what can we do when there is no geometric argument available to reduce the
problem to equations in one variable each?
In this case, we can still use a similar approach by replacing root-finding of
a polynomial in one variable by primary decomposition.

As our starting point, we take the $4$-parameter family $f_{s,t,u,v} \subset
\dP^3$ with dihedral symmetry $D_5$ defined by: 
\begin{equation}\label{eqnf}
\begin{array}{rcl}
p &:=& z \cdot \Pi_{j=0}^{4}\left[\cos\left({2\pi j\over
        7}\right)x 
    + \sin\left({2\pi j\over 7}\right)y -z\right]\\
& = & 
{z\over16}\Big[x\left(x^4-2\smcdot5\smcdot x^2y^2+5\smcdot y^4\right)\\
& & \qquad -5\smcdot z\smcdot\left(x^2+y^2\right)^2 
+4\smcdot5\smcdot z^3\smcdot\left(x^2+y^2\right)-16\smcdot z^5\Big],\\
q_{s,t,u,v} &:=& s\smcdot(x^2+y^2) +t\smcdot z^2 +u\smcdot zw + v\smcdot w^2,\\ 
f_{s,t,u,v} &:=& p-q_{s,t,u,v}^3.
\end{array}
\end{equation}

$p$ is the product of $z$ and $5$ planes in $\dP^3(\dC)$ meeting in the point 
$(0:0:0:1)$ with the symmetry $D_5$ of the $5$-gon with rotation axes
$\{x=y=0\}$.
$q_{s,t,u,v}$ is also $D_5$-symmetric, because $x$ and $y$ only appear as
$x^2+y^2$.   

The generic surface $f_{s,t,u,v}$ has $15\smcdot2=30$ singularities
of type $A_2$ at the intersections of the tripled quadric $q_{s,t,u,v}$ with
the ${6\choose2}$ pairwise intersection lines of the $6$ planes $p$. 
\ $2\smcdot5=10$ of the singularities lie in the $\{z=0\}$ plane, 
the other $4\smcdot5=20$ not. 
The coordinates of the latter $20$ can be obtained from the $4$ singularities
in the $\{y=0\}$ plane using the symmetry of the family. 
To see that the $\{y=0\}$ plane contains $4$ cusps, note that $p|_{y=0} =
z \cdot (z-x) \cdot (x^2-2xz-4z^2)^2$: For generic values of the parameters,
this doubled quadric factor meets the tripled quadric $q_{s,t,u,v}$ in
$2\cdot2$ points.  

Note that 
\begin{equation}\label{eqnv1}f_{s,t,u,v}(x,y,z,\lambda w) = f_{s,t,\lambda
  u,\lambda^2v}(x,y,z,w) \ \forall \lambda \in \dC^*,\end{equation}
s.t.\ we can choose
  $v:=1$ (it is easy to see that $v=0$ corresponds to a degenerate case).
Therefore, we write: 
$$f_{s,t,u} := f_{s,t,u,1} \ \textup{and} \ q_{s,t,u} := q_{s,t,u,1}.$$

To find surfaces in this $3$-parameter family with more singularities, we compute the
discriminant $Disc_{f_{s,t,u}}\in \dC[s,t,u]$ of the family $f_{s,t,u}$ by first 
dividing out the base locus (the intersections of the double lines of $p$
with the quadric $q$) from the singular locus (we use saturation, because we
have to divide out the base locus six times):
\begin{eqnarray*}
sl &:=& \left({\partial f \over \partial x}, {\partial f \over \partial y}, 
{\partial f \over \partial z}, {\partial f \over \partial w}\right),\\[0.5em]
bl &:=& \left({\partial p \over \partial x}, {\partial p \over \partial y}, 
{\partial p \over \partial z}, {\partial p \over \partial w},
q\right),\\[0.5em]
I &:=& sl \ : \ bl^\infty.
\end{eqnarray*}

Then we eliminate the variables $x,y,z$ from this quotient. 
In fact, because of the symmetry we restrict our attention to the
$\{y=0\}$ plane, which speeds up the computations: 
Every singularity in the plane $\{y=0\}$ which is not on the rotation axes
$\{x=y=0\}$ generates an orbit of length $5$ of singularities of the same
type. 
A short {\sc Singular} computation then gives the discriminant 
$Disc_{f_{s,t,u}} \in \dQ[s,t,u]$, which factorizes into $Disc_{f_{s,t,u}} =
D_{f,1} \cdot D_{f,2} \cdot D_{f,3}$, where: 
\begin{eqnarray*}
D_{f,1} &=& 2^{20} \smcdot 3^{6} \cdot  s^{5} \cdot \left(2^{4} \smcdot s^{2} + 2^{2} \smcdot 3 \smcdot st + t^{2}\right) \cdot \left(s + t\right)^{2}\\
&& + \left( - 2^{19} \smcdot 3^{6}\right) \cdot  s^{5} \cdot \left(2 \smcdot 11
 \smcdot s^{2} + 19 \smcdot st + 2 \smcdot t^{2}\right) \cdot \left(s +
 t\right) \cdot  u^{2}\\
&& +  2^{16} \smcdot 3^{6} \cdot  s^{5} \cdot \left(41 \smcdot s^{2} + 2 \smcdot 3 \smcdot 7 \smcdot st + 2 \smcdot 3 \smcdot t^{2}\right) \cdot  u^{4}\\
&& + \left( - 2^{14} \smcdot 3^{3}\right) \cdot  s^{3} \\
&& \quad \cdot \left(2 \smcdot
 3^{3} \smcdot 7 \smcdot s^{3}u^{6} + 2^{2} \smcdot 3^{3} \smcdot s^{2}tu^{6}
 + 2^{6} \smcdot 5^{2} \smcdot s^{3} - 2^{5} \smcdot 5^{2} \smcdot s^{2}t -
 5^{2} \smcdot 61 \smcdot st^{2} - 5^{3} \smcdot t^{3}\right) \\
& & +  2^{12} \smcdot 3^{3} \cdot  s^{3} \cdot \left(3^{3} \smcdot s^{2}u^{6}
 - 2^{5} \smcdot 5^{2} \smcdot s^{2} - 2 \smcdot 5^{2} \smcdot 61 \smcdot st -
 3 \smcdot 5^{3} \smcdot t^{2}\right) \cdot  u^{2} \\
&&  +  2^{10} \smcdot 3^{3} \smcdot 5^{2} \cdot  s^{3} \cdot \left(61 \smcdot s +
 3 \smcdot 5 \smcdot t\right) \cdot  u^{4} \\
&& + \left( - 2^{6} \smcdot 5^{3}\right) \cdot \left(2^{2} \smcdot 3^{3}
  \smcdot s^{3}u^{6} + 2^{6} \smcdot 5 \smcdot 23 \smcdot s^{3} + 2^{5}
  \smcdot 3 \smcdot 5 \smcdot s^{2}t + 2^{2} \smcdot 3 \smcdot 5^{2} \smcdot
  st^{2} + 5^{2} \smcdot t^{3}\right)\\ 
&& +  2^{4} \smcdot 3 \smcdot 5^{4} \cdot \left(2^{5} \smcdot s^{2} + 2^{3}
  \smcdot 5 \smcdot st + 5 \smcdot t^{2}\right) \cdot  u^{2}\\ 
&& + \left( - 2^{2} \smcdot 3 \smcdot 5^{5}\right) \cdot \left(2^{2} \smcdot s
 + t\right) \cdot  u^{4}\\ 
&& +  5^{5} \cdot \left(u^{4} - 2^{2} \smcdot u^{2} + 2^{4}\right) \cdot
 \left(u^{2} + 2^{2}\right),\\  
D_{f,2} &=& \left( - 2^{4}\right) \cdot  t^{2}
 +  2^{3} \cdot  t \cdot \left(u^{2} + 2\right)
 + \left(2 \smcdot u  - (u^{2} + 2^{2})\right) \cdot
 \left(2 \smcdot u  + (u^{2} + 2^{2})\right),\\ 
D_{f,3} &=& 2^{2} \cdot  t
 + \left(2 - u\right) \cdot \left(2 + u\right).
\end{eqnarray*}

We hope that some singularities of the discriminant correspond to examples of
surfaces $f_{s,t,u}$ with more $A_2$-singularities.
Note that only $D_{f,1}$ depends on the parameter $s$.  
Using computer algebra, it is easy to verify that the intersections of
two of the $3$ components $D_{f,1}, D_{f,2}, D_{f,3}$ of $Disc_{f_{s,t,u}}$ do
not yield to surfaces with many additional singularities.  

So, we use {\sc Singular} again to compute the primary decomposition
of the singular locus of $D_{f,1}$ over $\dQ$: $sl(D_{f,1}) =
\mathfrak{sl}_{f,1} \cap \mathfrak{sl}_{f,2} \cap 
\mathfrak{sl}_{f,3} \cap \mathfrak{sl}_{f,4}$ , where
$$\begin{array}{rclcl}
\mathfrak{sl}_{f,1} &=& \Big(2^{2} \smcdot \left(2^{2} \smcdot s - t\right) +
u^{2}, \quad 2^{6} \smcdot 3^{3} \smcdot  s^{3} - 5\Big)\\[0.4em]
\mathfrak{sl}_{f,2} &=& \Big( - 2^{2} \smcdot \left(2^{2} \smcdot 3
  \smcdot s + 5 \smcdot t\right) +  5 \smcdot  u^{2}, 
\quad 2^{4} \smcdot 3^{2} \smcdot  s^{2} +  2^{2} \smcdot 3 \smcdot 5 \smcdot  s +
5^{2} \Big)\\[0.4em]
\mathfrak{sl}_{f,3} &=& \Big( 2^{15} \smcdot 3^{3} \smcdot  t^{6}
 - 2^{14} \smcdot 3^{4} \smcdot  t^{5} \smcdot  u^{2}
 +  2^{11} \smcdot 3^{4} \smcdot 5 \smcdot  t^{4} \smcdot  u^{4}
 - 2^{6} \smcdot 3^{3} \smcdot 5 \smcdot  t^{3} \smcdot \left(2^{5} \smcdot
   u^{6} - 11 \smcdot 31\right)\\[0.2em] 
&&\qquad +  2^{4} \smcdot 3^{4} \smcdot 5 \smcdot  t^{2} \smcdot \left(2^{3} \smcdot
   u^{6} - 11 \smcdot 31\right) \smcdot  u^{2} 
 - 2^{2} \smcdot 3^{4} \smcdot  t \smcdot \left(2^{4} \smcdot u^{6} - 5
   \smcdot 11 \smcdot 31\right) \smcdot  u^{4} \\[0.2em]
&&\qquad  + \left(2^{3} \smcdot 3^{3} \smcdot u^{12} - 3^{3} \smcdot 5 \smcdot 11
   \smcdot 31 \smcdot u^{6} + 2^{6} \smcdot 5^{2} \smcdot 19^{3}\right),\\[0.4em]
&& \ \ 2^{11} \smcdot 3^{2} \cdot  t^{4}
 - 2^{11} \smcdot 3^{2} \cdot  t^{3} \cdot  u^{2}
 +  2^{8} \smcdot 3^{3} \cdot  t^{2} \cdot  u^{4}\\[0.2em]
&&\qquad  - 2^{2} \cdot \left(2^{5} \smcdot 3^{2} \smcdot tu^{6} - 2^{2} \smcdot 5
   \smcdot 7 \smcdot 19 \smcdot 211 \smcdot s - 5 \smcdot 73 \smcdot 193
   \smcdot t\right) \\[0.2em]
&&\qquad  + u^{2} \cdot \left(2^{3} \smcdot 3^{2} \smcdot u^{6} - 5 \smcdot 73 \smcdot
   193\right) 
\Big)\\[0.4em]
\mathfrak{sl}_{f,4} &=& \Big( 2^{2} \smcdot 3 \smcdot  s - 5, 
\quad -4 \smcdot (t+1) + u^2\Big).
\end{array}$$

All these prime ideals define smooth curves in the
$3$-dimensional parameter space. 
When projecting the curve $C_3$ defined by $\mathfrak{sl}_{f,3}$ to the $s,t$- or
the $s,u$-plane, we get in both cases six straight lines defined by the
equation \begin{equation}\label{eqns}2^{15} \smcdot 3^{3} \smcdot
  s^{6} - 2^{6} \smcdot 3^{3} \smcdot 5 \smcdot s^{3} + 5^{2}=0.\end{equation}
This shows that $C_3$ consists in fact of the union of six plane
curves. 
Over the algebraic extension $\dQ(s)$, it is easy to compute the equation of
these: 
\begin{equation}\label{eqnconics}C_{3,s}=5\smcdot u^2-2^2\smcdot5\smcdot t -
  2^{11} \smcdot 3^{2} \smcdot s^{4} - 2^{4} \smcdot 5 \smcdot s \in
  \dQ(s)[t,u].\end{equation} 

To show that there is a surface with $35$ $A_2$-singularities, we take the
most simple point of this curve, the one with $u=0$:
\begin{theorem}\label{thm35}
Let $s_0\in\dC$ be one of the six roots of (\ref{eqns}).
Let $(t_0, 0)$ be the point on $C_{3,s_0}$ with $u=0$.    
Then the sextic $f_{s_0,t_0,0} \subset \dP^3$ has exactly $35$
singularities of type $A_2$ and no other singularities.
\end{theorem}
\begin{proof}
We use computer algebra. 
The {\sc Singular} script 
and its output 
can be downloaded from the webpage
\cite{labsAlgSurf}.  
Here, we give the basic ideas:

%
With $u=0$ in $C_{3,s_0}$, we find: $t_0 =
-4\smcdot s_0\left({2^{7} \smcdot 3^{2} \over 5} \smcdot s_0^{3} + 1\right).$
For the corresponding surface 
\begin{equation}\label{eqnS35}
S_{35}:= f_{s_0, -4\smcdot s_0\left({2^{7} \smcdot 3^{2} \over 5} \smcdot s_0^{3} +
    1\right), 0}\end{equation}
we first check that the total milnor number is $70$. 
Then we verify that the surface has $35$ singularities of type $A_2$: 
For each orbit of singularities, we compute the ideal of one of the
singularities and check explicitly that it is a cusp.
To show this it suffices to verify that its milnor number is exactly two. 
E.g., for the orbit of the five non-generic singularities, we take the cusp
    $S_{yw}$ that lies in the $\{y=0\}$ plane with coordinates:  
$$S_{yw} = \left(-{2^7\smcdot 3^2\over5}s_0^3 +8 : 0 : 1 : 0\right).$$
\end{proof}

The previously known lower bound for the number of $A_2$ singularities a
sextic surface in $\dP^3$ can have, was $30$ (using Rohn's construction
mentioned in the introduction).
We now have, with Miyaoka's upper bound:
\begin{corollary}
The maximum number $\mu_{A_2}(6)$ of cusps a surface of degree
$6$ in $\dP^3$ can have, satisfies:
$$35 \le \mu_{A_2}(6) \le 37.$$
\end{corollary}



Note that the coefficients of the surface $S_{35}$ are not real. 
In fact, the ideal $\mathfrak{sl}_{f,3}$ does not contain any real
point, because equation (\ref{eqns}) does not have any real root.
In particular, it is not possible to use the software {\tt surf} \cite{surf}
to draw an image of this sextic.
This also holds for the more general family $f_{s,t,u,v}$ because of equation
(\ref{eqnv1}). 
The curves defined by the ideals $\mathfrak{sl}_{f,2}$ and
$\mathfrak{sl}_{f,4}$ lead to only one additional higher singularity, and we
are not interested in such examples.  

\begin{figure}[htbp]
\begin{center}
\includegraphics[width=1.5in]{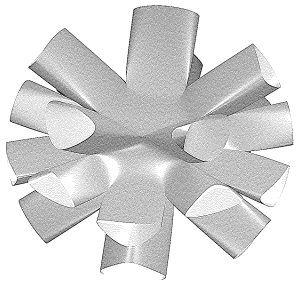}
\end{center}
\caption{A sextic with $30$ cusps and $10$ nodes at infinity.}
\label{fig3010}
\end{figure}

But in the case of the prime ideal $\mathfrak{sl}_{f,1}$, 
we get surfaces with $30$ real $A_2$-singularities and $10$ real
$A_1$-singularities (see also fig.~\vref{fig3010}). 
Again, we choose a point in the parameter-space with $u=0$:
\begin{theorem}
The sextic $f_{s_0,t_0,0} \subset \dP^3$, where $s_0 :=
{1\over3\smcdot2^2}\sqrt[3]{5}\in\dR$, $t_0 = 2^2\smcdot s_0\in\dR$, has exactly $30$
singularities of type $A_2$, $10$ 
singularities of 
type $A_1$, and no other singularities.  
Furthermore, all the singularities are real. 
\end{theorem}
\begin{proof}Similar to the preceding one.\end{proof}



\section{Concluding Remarks}

The method we used to find the sextics with $35$ cusps within a family of
surfaces with many singularities is an algorithm that can be applied to many
other families of hypersurfaces.  
The only limit is the time needed for the computations. 
The idea to such an algorithm is not new. 
In fact, our main observation was to notice that we can use features of the
most recent versions of the computer algebra system {\sc Singular}
to perform the algorithm on a computer in our particular
case:
Finding the equation of $S_{35}$ and verifying that it has exactly $35$ cusps
and no other singularities just takes a few seconds on our computer. 

We could restrict our attention to a plane because of the
symmetry of our family, so that the number of variables decreased. 
This speeded up the computations. 
But the case of septics with many nodes was too time-consuming to be
treated in this way:
Our construction of a $99$-nodal surface of degree $7$ involves 
computations in positive characteristics and then liftings to characteristic
zero using the geometry of the examples, see \cite{labs99}. 

In other applications, it might be easy to divide out the base locus and to
compute the discriminant, e.g.\ by using the geometry of the
family. 
Then it only remaines to study the discriminant for
finding examples which have more singularities than the generic member of the
family.


\nocite{gpSingComAlg}
\nocite{labsAlgSurf}

\bibliographystyle{plain} 
\bibliography{papers}

\end{document}